\documentclass[12pt]{article}
\usepackage{amsmath,amssymb,amsbsy,amsfonts,amsthm,latexsym,
               amsopn,amstext,amsxtra,euscript,amscd}
\usepackage[francais]{babel}
\usepackage[latin1]{inputenc}

\newcommand\Retire[1]{\relax}

\begin{document}

\newtheorem{theorem}{Theorem}
\newtheorem{lemma}[theorem]{Lemma}
\newtheorem{claim}[theorem]{Claim}
\newtheorem{cor}[theorem]{Corollaire}
\newtheorem{proposition}[theorem]{Proposition}
\newtheorem{definition}{Définition}
\newtheorem{question}[theorem]{Open Question}
\newtheorem{ques}[theorem]{Question}
\newtheorem{theo}{Th\'eor\`eme}
\newtheorem{lemme}[theorem]{Lemme}
\def\E{{\mathbf E}}
\def\cA{{\mathcal A}}
\def\cB{{\mathcal B}}
\def\cC{{\mathcal C}}
\def\cD{{\mathcal D}}
\def\cE{{\mathcal E}}
\def\cF{{\mathcal F}}
\def\cG{{\mathcal G}}
\def\cH{{\mathcal H}}
\def\cI{{\mathcal I}}
\def\cJ{{\mathcal J}}
\def\cK{{\mathcal K}}
\def\cL{{\mathcal L}}
\def\cM{{\mathcal M}}
\def\cN{{\mathcal N}}
\def\cO{{\mathcal O}}
\def\cP{{\mathcal P}}
\def\cQ{{\mathcal Q}}
\def\cR{{\mathcal R}}
\def\cS{{\mathcal S}}
\def\cT{{\mathcal T}}
\def\cU{{\mathcal U}}
\def\cV{{\mathcal V}}
\def\cW{{\mathcal W}}
\def\cX{{\mathcal X}}
\def\cY{{\mathcal Y}}
\def\cZ{{\mathcal Z}}
\def\E{{\mathbb E}}
\def\D{{\mathbb D}}
\def\F{{\mathbb F}}
\def\Z{{\mathbb Z}}
\def\N{{\mathbb N}}
\def\Q{{\mathbb Q}}
\def\R{{\mathbb R}}
\def\F{{\mathbb F}}
\def\K{{\mathbb K}}
\def\L{{\mathbb L}}
\def\M{{\mathbb M}}
\def\Z{{\mathbb Z}}
\def\A{{\mathbb A}}
\def\B{{\mathbb B}}
\def\U{{\mathbb U}}
\def\e{{\mathbf e}}
\def\a{{\bf a}}
\def\b{{\bf b}}
\def\cc{{\bf c}}
\def\d{{\bf d}}
\def\h{{\bf h}}
\def\ee{{\rm e}}

\def\\{\cr}
\def\({\left(}
\def\){\right)}
\def\lcm{{\rm lcm\/}}
\def\fl#1{\left\lfloor#1\right\rfloor}
\def\rf#1{\left\lceil#1\right\rceil}
\def\ov{\overline}

\title{Sur la complexité de familles d'ensembles pseudo-aléatoires}
\author{
{\sc R. Balasubramanian (Chennai)}\and {\sc C\'ecile Dartyge (Nancy)}\and {\sc \'Elie Mosaki (Lyon)}}

\date{\today}

\maketitle

\renewcommand{\thefootnote}{}
\footnote{\noindent\textbf{} \vskip0pt \textbf{Mots clés :} sous-ensembles pseudo-aléatoires, complexité, sommes d'exponentielles, sommes de caractères.

\noindent\textbf{} \vskip0pt \textbf{Classification . :} 11K45, 11L07, 05B10}

\begin{abstract}
Dans cet article on s'intéresse au problème suivant.
Soient $p$ un nombre premier, $S\subset \F_p$ et $\cP\subset \{ P\in\F_p [X] : \deg P\le d\}$.
Quel est le plus grand entier $k$ tel que pour toutes paires de  sous-ensembles disjoints $\cA ,\cB$ de $\F_p$ vérifiant 
$|\cA\cup\cB |=k$, il existe $P\in\cP$ tel que $P(x)\in S$ si $x\in\cA$ et $P(x)\not\in S$ si $x\in\cB$?
Ce problème correspond à l'étude de la complexité de certaines familles d'ensembles pseudo-aléatoires.
Dans un premier temps nous rappelons la définition de cette complexité et resituons le contexte
 des ensembles pseudo-aléatoires.
 Ensuite nous exposons les différents résultats 
obtenus selon la nature des ensembles $S$ et $\cP$  étudiés.
Certaines preuves passent par des majorations de sommes d'exponentielles ou de caractères sur des corps finis,
d'autres combinent des arguments combinatoires avec des résultats de la théorie additive des nombres.
\medskip

\centerline{\bf Abstract}
In this paper we are interested in the following problem. Let $p$ be a prime number, $S\subset \F_p$ and $\cP\subset \{ P\in\F_p [X] : \deg P\le d\}$. What is the largest integer $k$ such that for all subsets  $\cA, \cB$ of $\F_p$ satisfying $\cA\cap\cB =\emptyset$ and
$|\cA\cup\cB |=k$, there exists  $P\in\cP$ such that $P(x)\in S$ if $x\in\cA$ and $P(x)\not\in S$ if $x\in\cB$?
This problem corresponds to the study of the complexity of some families of pseudo-random subsets.
First we recall this complexity definition and the context of pseudo-random subsets. Then we state the different results 
we have obtained 
according to the shape of the sets $S$ and $\cP$ considered.
 Some proofs are based on upper bounds for exponential sums or characters sums in finite fields, other proofs 
 use combinatorics and additive number theory.
\end{abstract}

\section{Introduction}

Dans des problèmes de simulation ou de cryptographie, on a parfois besoin de sous-ensembles de $\{ 1,\ldots , N\}$ ou de $\Z _n$ qui ressemblent à  des ensembles d'entiers pris au hasard.
Dartyge et  S\'ark\"ozy   \cite{DS07a} et \cite{DS07b} ont
proposé pour cela,  des mesures de nature pseudo-aléatoire de sous-ensembles de
$\{ 1,\ldots ,N\}$ et de $\Z _n$, où $n$ et $N$ sont des entiers donnés.
Ces mesures reprennent celles de bonne corrélation et de bonne répartition dans les progressions arithmétiques
définies par Hubert, Mauduit et S\'ark\" ozy \cite{MS97} et \cite{HS04} pour les suites binaires pseudo-aléatoires.

Soit $\cR\subset \{ 1,\ldots ,N\}$. On associe à  $\cR$ le $N$-uplet $(e_1,\ldots,e_N)$~:
\begin{displaymath}
e_n=\left\{\begin{array}{ll} 1-\frac{|\cR |}{N} & {\text si}\ n\in\cR \\
-\frac{ |\cR|}{N} & {\text si}\ n\not\in\cR
\end{array}\right. \quad (n=1,\ldots , N).
\end{displaymath}
La mesure de bonne répartition dans les progressions arithmétiques est alors
$$W(\cR ,N)=\max _{a,b,t}\Big |\sum_{j=0}^{t-1}e_{aj+b}\Big |,$$
où le maximum porte sur les entiers $a,b,t$ tels que $1\le b\le b+(t-1)a\le N$. La seconde mesure est la mesure de corrélation d'ordre $k$ :
$$C_k(\cR , N)=\max_{M,D}\Big |\sum_{n=1}^M e_{n+d_1}\cdots e_{n+d_k}\Big |,$$
où le maximum est sur les $D=(d_1,\ldots ,d_k)$ et $M$ tels que
$0\le d_1<\cdots <d_k\le N-M$. Dans le cas des sous-ensembles de $\Z_n$ les conditions sur $a,b,D$ pour ces deux mesures sont légèrement différentes.

On dira alors que $\cR$ possède de bonnes propriétés pseudo-aléatoires, si
$W(\cR ,N)=o(N)$ et $C_k(\cR , N)=o(N)$.

S\'ark\"ozy, Szalay et  les deux derniers auteurs ont donné dans les articles
\cite{DMS}, \cite{DS07b} et \cite{DSS09} plusieurs exemples de familles d'ensembles pseudo-aléatoires.
Dans certaines applications il est  important de disposer
de familles d'ensembles aléatoires dont la structure est riche, ou encore dont les éléments ne peuvent pas  être déterminés à  partir d'un faible nombre de données.
Cela a conduit  Ahlswede,  Khachatrian, Mauduit et
 S\'ark\"ozy \cite{AKMS03} à  définir la notion de complexité d'un ensemble de suites pseudo-aléatoires.

Dans les articles \cite{DMS} et \cite{DS07b}
cette définition a été adaptée dans le cadre
de  familles de sous-ensembles pseudo-aléatoires de $\{ 1,\ldots ,N\}$
de la manière suivante~:
\begin{definition}
\label{def1}
Soit  $\cF$ une famille   de sous-ensembles de  $\{1,2,\ldots,N\}$.
La  complexité $K (\cF )$ de la famille $\cF$  est le plus grand entier
 $k\in \N$ tel que pour tout  ${\cal A}\subset
\{1,2,\ldots,N\}$ avec $|{\cal A}|=k$ et tout sous-ensemble ${\cal B}$
de ${\cal A}$  il existe un élément de ${\cal F}$ tel que ${\cal R}\cap{\cal A}={\cal
B}$. Autrement dit, pour tout  ${\cal A}\subset \{1,2,\ldots,N\}$
tel que $|{\cal A}|=k$ et toute  partition
$${\cal A}={\cal B}\cup{\cal C}\;,\;{\cal B}\cap{\cal C}=\varnothing$$
de ${\cal A}$ il existe   ${\cal R}\in {\cal F}$ tel que
$${\cal B}\subset{\cal R},\, {\text et}\
{\cal C}\subset \{1,\ldots,N\}\setminus{\cal
R}\,.$$
\end{definition}

L'objet de  cet article est de  continuer l'étude de la complexité de différentes familles étudiées dans \cite{DMS} et \cite{DS07b}.

Les ensembles construits dans \cite{DMS} et \cite{DS07b} sont de la forme :
$$\cR (f,S)=\{ n\in\{ 1,\ldots ,p\} : \exists h\in S\ {\rm tel}\ {\rm que}\ f(n)\equiv h \mod p\},$$
où $f$ est un polynôme et $S$ est un ensemble non vide strictement contenu dans $\F_p$.

Les ensembles  $S$ considérés   étaient des différents types suivants :

(i) $S_1=\{ r,r+1,\ldots ,r+s-1\}$ avec $r\in\F_p$ et $s<p/2$ (\cite{DMS}) ;

(ii) $S_2=\{ \bar r,\ov{r+1},\ldots ,\ov{r+s-1}\}$ avec la notation
$x\ov x\equiv 1\mod p$  (\cite{DMS}), en omettant $\ov{0}$ par convention;

(iii) $S_3$ est l'ensemble des puissances $\ell $ ièmes modulo $p$ où $\ell$ est un diviseur de $p-1$ (\cite{DS07b}) .

Les ensembles de type $\cR ( f, S_1)$ ont l'avantage  d'être rapides
à  générer mais peuvent avoir des mauvaises mesures de corrélation (cf. \cite{MRS04} Théorème 4).

Pour les ensembles de type $S_2$ ou $S_3$, les mesures de corrélation sont plus
difficiles à   estimer. Dans \cite{DMS} et \cite{DS07b} ces mesures ont été évaluées pour des polynômes de la forme $f_\cA (X)=\prod_{a\in\cA}(X-a)$ où
$\cA \subset \F _p$ ou plus généralement pour des polynômes $f$ sans racine multiple mais avec en contrepartie des conditions plus contraignantes
sur les rapports des  corrélations.

Notons $\cP_1 (d,p)$ l'ensemble des polynômes de $\F _p[X]$ de degré $\le d$,
$\cP_2 (d,p)$ celui des polynômes sans racine multiple et de degré   $\le d$
et enfin $\cP _3 (d,p)=\{ f_\cA :\cA\subset\F_p, |\cA |= d\}$ avec la notation $f_A$
définie ci-dessus. On a ainsi les inclusions $\cP_3(d,p)\subset\cP_2(d,p)\subset \cP_1(d,p)$.

On définit également pour i=1,2,3 :
$$\cF _i(S,d)=\{ \cR (f, S) : f\in\cP_i (d,p)\},$$
et $K_i (S,d)$ la complexité correspondante.

En utilisant les polynômes d'interpolation de Lagrange, on voit facilement que
$K_1 (S,d)\ge d+1$.
Dans \cite{DMS} nous montrons à  l'aide du théorème de Cauchy-Davenport
que $K_1 (S,d)\ge d+2$.
Ce dernier résultat est valable pour tous les ensembles $S$ tels que $\min (|S|, |S^c|)$ est assez grand ;
$|S^c|$ étant le complémentaire de $S$ dans $\F_p$ : 
$S^c=\F_p\setminus S$. On pourrait penser que
dans le cas où $S$ est une suite d'entiers consécutifs l'on puisse obtenir
une meilleure minoration. Nous n'y sommes pas parvenus.
Les sommes d'exponentielles associées à  ces problèmes se calculent de manière élémentaire et sont dans certains cas de taille très importante.
Ce phénomène apparaît également  dans l'étude des mesures de corrélations des ensembles $\cR (f,S_1)$ (cf. \cite{MRS04}).

Par contre, la situation est différente
lorsque
$S$ est l'ensemble des inverses d'une suite d'entiers consécutifs.

 \begin{theo}
\label{thm:K3} Soient $d\ge 2$, $k\in\N^*$, $r\in\F_p$ et $\beta\in ]0,1[$ donnés.
On note $s=\lceil \beta p\rceil$.
On considère le sous-ensemble $S\subset\F_p$ défini par
$$S=\{ \ov r,\ov{r+1},\ldots ,\ov{r+s-1}\}\quad (on\; omet\; \ov{0}).$$
L'inégalité
$K_3(S,d)\ge k$ est alors vérifiée pour  $p>2/(1-\beta)$ tel que
\begin{equation}\label{condK3p}
{p-k\choose d}\min_{0\le \ell\le k}\beta^\ell\Big (1-\beta-\frac{1}{p}\Big )^{k-\ell}-
(44^{4k}+2k+2d-2)p^{d-1}(\log (9p))^k>0.
\end{equation}
\end{theo}
Comme ${p\choose d}\sim \frac{p^d}{d!}$, on en déduit que $K_3(S,d)\ge k$
pour $p$ assez grand, $p\ge p_0(\beta ,d,k)$.
Il découle de ce théorème la minoration :
\begin{equation}
\label{minK3}
K_3(S,d)\gg_{d,\beta}\frac{\log p}{\log\log p}.
\end{equation}
Rappelons que  la proposition 4.3 de \cite{DMS} entraîne que $K_3(S,d)\le \frac{(d+1)\log p}{\log 2}$.

Plus $d$ est grand  plus la condition \eqref{condK3p} est mauvaise. Cela est contraire à  notre intuition. Ce défaut est dû aux majorations de sommes d'exponentielles que nous utilisons.

Pour $y\in\R$ et $p$ premier on note $\e _p(x)=\exp (2i\pi x/p)$.
L'une des étapes de la preuve de ce théorème est de trouver des majorations de sommes de la forme :
$$S:=\sum_{\{a_1,\ldots ,a_d\}\subset\F_p}\e_p \Big ( \sum_{m=1}^th_m
\prod_{j=1}^d\ov{b_{m}-a_j}\Big ),$$
où chaque $a_j$ est  différent des $b_{m}$.
Eichenauer-Herrmann et  Niederreiter \cite{E-HN94} ont étudié ces sommes dans le cas $d=1$. En utilisant les majorations de Bombieri et Weil  (\cite{Bo66})
de sommes d'exponentielles d'argument une fraction rationnelle,
ils obtiennent des majorations  complètement explicites.

On déduit facilement du  théorème d'Eichenauer-Herrmann et Niederreiter la majoration $S\ll_{\beta,k,d} p^{d-1/2}$, où la constante implicite est calculable. Cette majoration est suffisante pour obtenir
une version légèrement affaiblie du
Théorème \ref{thm:K3}. Le deuxième terme du membre de gauche de \eqref{condK3p} étant alors de l'ordre de $p^{d-1/2}(\log p)^k$.

Une question naturelle est  de vérifier s'il n'est pas possible d'obtenir une meilleure majoration de $S$ en profitant du fait d'avoir une somme sur $d$ variables avec $d\ge 2$. Lorsque $t=1$, nous verrons au paragraphe 2 que la somme $S$ s'évalue facilement et est de l'ordre de $p^{d-1}$. Nous  pouvons donc profiter de compensations seulement sur deux variables.

Le théorème suivant est  on l'espère d'un intérêt intrinsèque.
Il  est un résultat  analogue au théorème d'Eichenauer-Herrmann et Niederreiter pour des sommes en deux variables.
\begin{theo}
\label{expo}
Soient $b_1,\ldots , b_k$, $c_1,\ldots ,c_k$ des éléments de $\F_p$
tels que $b_i\not =b_j$ et $c_i \not = c_j$ pour tous $i\not =j$.
Pour $d_1,\ldots ,d_k\in\F_p^*$,
on considère la somme d'exponentielles :
$$S(\b,\cc, \d)=\sum_{\substack{x,y\in\F_p\\ x\not =b_i, y\not =c_i}}
\e_p(\sum_{i=1}^kd_i\ov{(x-b_i)(y-c_i)}).$$
On a alors
\begin{equation}
\label{majSbcd}
|S(\b,\cc,\d)|\le 44^{4k}p.
\end{equation}
\end{theo}
Une application directe du théorème d'Eichenauer-Herrmann et Niederreiter
fournit la majoration $S(\b ,\cc,\d)=O(p^{3/2}).$
Nous verrons au paragraphe 2 que le Théorème \ref{expo}
implique la majoration $S\ll p^{d-1}$.
Il faut cependant signaler que la dépendance en $k$ de notre théorème est de moins bonne qualité que celle 
issue de la majoration d'Eichenauer-Herrmann et Niederreiter. Autrement dit, lorsque $k$ est grand comparativement à $p$
(en fait pour $k\gg \log p$)
il est plus pertinent d'utiliser le Théorème d'Eichenauer-Herrmann et Niederreiter.

La preuve de ce résultat fait l'objet du premier paragraphe de cet article.
Elle utilise des résultats de géométrie algébrique notamment les célèbres travaux de Dwork et de Deligne. Nous avons essayé d'adopter une approche
la plus élémentaire possible en nous inspirant des travaux de Hooley \cite{H82}, Adolpshon et Sperber \cite{AS89} ainsi que de la présentation
de la méthode de Hooley faite par Birch et Bombieri \cite{BiBo85}
dans un appendice d'un article de Friedlander et Iwaniec \cite{FI85}.

Cette preuve nécessite l'étude d'extensions des sommes $S(\b ,\cc,\d)$ à des sommes d'exponentielles définies sur $\F_{p^n}$
et fournit au passage des majorations de telles sommes.

\medskip
La structure de $S$ est essentielle dans la preuve du Théorème \ref{thm:K3}.
Pour un ensemble $S$ général c'est plus délicat. On ne peut même pas utiliser les polynômes d'interpolation de Lagrange car
rien ne dit que le polynôme alors formé sera dans $\cP_3(d,p)$ (c'est-à-dire à racines simples dans $\F_p$).
Le résultat suivant donne une minoration de $K_3(S,d)$ valable pour des ensembles $S$ quelconques.

\begin{theo}
\label{k3general}
(i) Si $S$ et $S^c$ sont non vides alors
$K _3(S,d)> \lceil d/2 \rceil$ pour $p$ assez grand ($p>p_0 (d)$).
\par
(ii) si $\min (|S|, |S^c|)\gg p$ alors $K_3(S,d)\ge d-1$ pour $p$ assez grand.
\end{theo}

Pour démontrer ce théorème nous utilisons des majorations de sommes de caractères multiplicatifs
 afin de profiter de la structure produit des polynômes $f_\cA$.

La fin de cet article est dévolue à  la complexité $K_2(S,d)$.
Dans \cite{DS07b} et \cite{DMS} on trouve deux arguments différents montrant que cette complexité est $K_2 (S,d)\ge d+1$ si $S$ et son complémentaire ont suffisamment d'éléments.
Dans cet article nous améliorons ce résultat :

\begin{theo}
\label{k2general}
Si $4d-2<|S|<\frac{p-d+1}{2}$ alors
$K_2 (S,d)\ge d+2$.

\end{theo}

La preuve de ce résultat est de nature différente de celle des précédents théorèmes. Elle est combinatoire et utilise un résultat récent de Green et Ruzsa \cite{GR} de théorie additive des nombres.
Pour $x\in\F _p$, on note $\bar x$ l'inverse de $x$ dans $\F_p$. Pour $n\in\N$, (resp. $n\in\F_p)$, on note $r_p(n)$ le plus petit entier positif congru à  $n$ modulo $p$ (resp. appartenant à  la classe de $n$ dans $\F _p$).
En fait dans cet article nous ferons souvent l'amalgame entre un entier $n$ et sa classe dans $\F_p$.

\section{Majorations de sommes d'exponentielles}

Soient $\h\in (\F_p^*)^k$, $\b =(b_{m,j})_{\substack {1\le m\le k\\
1\le j\le d}}\in\cB^{kd}$ où $\cB\subset\F_p$.

Dans ce paragraphe nous étudions des sommes d'exponentielles de la forme :
\begin{equation}\label{Shb}S(\h,\b)=\sum_{ \{ a_1,\ldots ,a_d\}\subset\F_p\setminus \cB}
\e _p \Big (
\sum_{m=1}^{k}h_m\prod_{j=1}^d \ov{ b_{m,j}-a_j}
\Big ),\end{equation}
où pour chaque $j$, $b_{m,j}\not = b_{\ell,j}$ si $\ell\not =m$.
\subsection{Sommes sur une variable}
Pour $q=p^n$, $x\in\F_q$, on note encore $\ov x$ l'inverse de $x$ dans $\F_q$.
Commençons par rappeler le résultat d'Eichenauer-Herrmann et Niederreiter :

\begin{theo}
\label{d1} {\bf (Eichenauer-Herrmann et Niederreiter \cite{E-HN94} Théorème 1 p. 270.)}
Soient $\d\in\F_q^s$, $\d\not =0$ et $\e =(e_1,\ldots ,e_s)\in \F_q^s$
tels que $e_1,\ldots ,e_s$ soient distincts deux à  deux. Si $\psi$ est un caractère additif de $\F_q$ non trivial alors
$$\sum_{n\in\F_q\setminus\{ -e_1,\ldots ,-e_s\}}\psi\Big (\sum_{j=1}^s d_j
\ov{n+e_j}\Big )\le (2s-2)\sqrt{q} +1.$$
\end{theo}
Cette version est légèrement différente du théorème 1 de \cite{E-HN94}
car dans cet article la somme porte sur tous les $n\in\F_q$ avec la convention $\ov x=0$ pour $x=0$ ce qui crée un $s$ en plus dans la majoration.

Dans le cas où $k=1$, et $b_{1,j}=b$ pour $1\le j\le d$, la somme $S(\h ,\b)$ définie par
\eqref{Shb} est simplement du type :
$$S(h,b)=\sum_{ \{ a_1,\ldots ,a_d\}\subset\F_p\setminus \{ b\}}
\e_p(h\prod_{j=1}^d\ov{b-a_j}).$$

\begin{lemme}
\label{t1}
Soit $d\ge 2$. On a l'égalité
$$S(h,b)=-\frac{p^{d-1}}{(d-1)!}+O_d(p^{d-3/2}).$$

\end{lemme}
La preuve de ce lemme repose sur le fait que pour $(a,p)=1$,
\begin{equation}
\label{orthcarac}\sum_{x\in\F_p^*}\e_p (a\ov x)=\sum_{x\in\F_p}\e_p(ax)-1=-1.
\end{equation}
On commence par sommer sur $a_d\in\F_p\setminus
\{b,a_1,\ldots ,a_{d-1}\}$, $a_1,\ldots ,a_{d-1}$ étant fixés.
En complétant la somme sur $a_d$  pour utiliser \eqref{orthcarac} on obtient :

\begin{displaymath}
S(h,b)=-{p-1\choose d-1}-\sum_{i=1}^{d-1}\sum_{\{a_1,\ldots ,a_{d-1}\}\subset\F_p\setminus \{b\}} \e_p(h\ov{(b-a_i)^2}\prod_{\substack{j=1\\ j\not =i}}^{d-1}\ov{b-a_j}).
\end{displaymath}
Le premier terme est  de l'ordre de $p^{d-1}/(d-1)!$ tandis que les autres termes de la somme en $i$, qui sont tous égaux, sont des $O(p^{d-3/2})$, d'après le Théorème \ref{d1} pour $d\ge 3$
ou d'après les majorations classiques de sommes de Gauss pour $d=2$.
Cela termine la preuve du Lemme \ref{t1}.

\medskip

{\bf Remarque.} Lorsque $d\ge 3$, on peut améliorer le terme d'erreur du lemme  \ref{t1} en appliquant
\eqref{orthcarac} un nombre approprié de fois.

\subsection{Premières étapes de la preuve du Théorème \ref{expo}.}

Si $k=p$ le résultat est évident. 

\medskip 

Soit $k<p$. Quitte à  faire des changements de variables on peut supposer
\begin{equation}
\label{hypbc}
\prod_{i=1}^kb_ic_i\not= 0.
\end{equation}
En effet supposons que cette condition ne soit pas réalisée.
Soient $\beta\in\F_p\setminus\{ b_1,\ldots ,b_k\}$ et $\gamma\in\F_p\setminus\{ c_1,\ldots ,c_k\}$.
En posant $x'=x-\beta$, $y'=y-\gamma$, la somme devient :
$$S(\b ,\cc,\d)=\sum_{\substack{x'\in\F_p\setminus\{b_1-\beta,\ldots ,
b_k-\beta\}\\ y'\in\F_p\setminus\{ c_1-\gamma,\ldots , c_k-\gamma\}}}
\e_p\Big (\sum_{i=1}^k d_i\ov{(x'+\beta -b_i)(y'+\gamma -c_i)}\Big ),$$
où maintenant  $b_i-\beta\not =0$, $c_i-\gamma\not =0$ pour $1\le i\le k$ ;
on s'est ramené à  une somme vérifiant \eqref{hypbc}.

\medskip

Posons $u_i=\ov{x-b_i}$, $v_i=\ov{y-c_i}$ pour $1\le i\le k$.
On a alors pour $1\le i\le k$ :
\begin{equation}
\label{eqV}
u_iu_1(b_1-b_i)+u_i-u_1=0\ {\rm et}\ v_iv_1(c_1-c_i)+v_i-v_1=0.
\end{equation}
Soit $V\subset (\F_p)^{2k}$ la variété définie par les équations \eqref{eqV}.
\'Etant donné que tous les $u_i$ s'expriment en fonction de $u_1$
et tous les $v_i$ en fonction de $v_1$, $V$ est une variété de dimension 2.

La somme d'exponentielles $S(\b,\cc,\d)$ se réécrit alors  de la manière suivante :
\begin{equation}
\label{SV}S(\b ,\cc ,\d)=\sum_{\substack{u,v\in(\F_p^*)^k\\ (u,v)\in V}}
\e_p\big (\sum_{i=1}^kd_iu_iv_i\big ).
\end{equation}
La somme $S(\b,\cc,\d)$ apparaît maintenant comme une somme d'exponentielles de la forme :
$$\sum_{\substack{0\le x_1,\ldots ,x_n<p\\ g_1(x_1,\ldots ,x_n)=\cdots
=g_s (x_1,\ldots ,x_n)=0\ {\rm mod}\ p}}\e_p (f(x_1,\ldots ,x_n)),$$
avec $f,g_1,\ldots ,g_s\in \F_p[X_1,\ldots,X_n]$.
L'étude de ce type de  sommes est une branche  de la  géométrie arithmétique qui  connaît un développement très important depuis le siècle dernier.
Lorsque la somme ne porte que sur une variable, on dispose de majorations valables dans un cadre très général grâce aux travaux de Weil, puis de Deligne et Bombieri.
Lorsque la somme porte sur plusieurs variables la situation est moins connue.
Dans le cas où la somme est de la forme $\sum_{x_1,\ldots ,x_n\in\F _p}
\e_p( f(x_1,\ldots ,x_n))$, où $f$ est de degré $d$
et a une composante homogène de degré $d$ non singulière, Deligne \cite{De74} a montré que cette somme est de module inférieur à 
$(d-1)p^{n/2}$.

D'autres résultats très profonds ont été obtenus ces dernières décennies.
Rojas-Le\'on \cite{R-L07} a par exemple récemment établi des majorations
de telles sommes avec des hypothèses sur $f$ moins fortes.
On trouvera dans cet article d'autres références sur cette question.
Malheureusement nous n'avons pas pu appliquer le résultat de Rojas-Le\'on évoqué ci-dessus.

Hooley \cite{H82} a repris les travaux de Dwork et de Deligne et a
proposé une méthode pour obtenir des majorations dans un cadre assez général.
Nous nous sommes inspirés de son approche ainsi que de la présentation qui en est faite
dans l'appendice de Birch et Bombieri \cite{BiBo85} d'un article de Friedlander et Iwaniec \cite{FI85}.

Pour $n\in \N$, on note $\F_{p^n}$ une extension de $\F_p$ de dimension $n$.
Pour tout $x\in\F_{p^n}$ on note $\sigma _n(x)$ la trace de $x$ sur $\F_p$ :
$$\sigma_n(x)=x+x^p+\cdots +x^{p^{n-1}}.$$
Rappelons que  $\sigma_n(u)\in\F_p$ si $u\in\F_{p^n}$.

Soit $V_n$ l'ensemble des $(x,y)\in\F_{p^n}^{2k}$ vérifiant \eqref{eqV}.
On forme ensuite les sommes d'exponentielles sur $\F_{p^n}$ et $\F_{p^n}^*$:
$$S_n(V, \d)=\sum_{(x,y)\in V_n}\e_p \big (\sigma _n
\big (\sum_{i=1}^kd_ix_iy_i\big )\big ),$$
et
$$S_n^*(V, \d)=\sum_{\substack{(x,y)\in V_n\\ x,y\in(\F_{p^n}^*)^k}}\e_p \big (\sigma _n
\big (\sum_{i=1}^kd_ix_iy_i\big )\big ).$$

Si $(x,y)\in V_n$ est tel que $x_i=0$ pour un $i\in\{ 1,\ldots ,k\}$ donné
alors d'après les équations de $V$, tous les $x_j$ sont nuls.
De même si l'un des $y_j$ est nul tous les autres le sont.
On en déduit l'égalité
\begin{equation}
\label{SnSn*}S_n(V,\d)=S_n^*(V,\d)+2(p^n-1)+1.
\end{equation}

On considère alors les séries de Dirichlet  définies formellement par
$$L(T)=\exp\Big ( \sum_{r=1}^\infty \frac{S_r(V,\d)T^r}{r}\Big )\ {\rm et}\
L^*(T)=\exp\Big ( \sum_{r=1}^\infty \frac{S_r^*(V,\d)T^r}{r}\Big ).$$
Dwork \cite{Dw60} et Bombieri \cite{Bo66} ont montré
que les séries
$L(T)$ et $L^*(T)$ sont des fractions rationnelles dont les numérateurs et les dénominateurs appartiennent à  $\Q(\e_p(1))[T]$.
Cette preuve est reprise dans l'article de Hooley \cite{H82}.
Il en déduit ensuite pour chaque entier $r$ l'égalité :
\begin{equation}
\label{SrVd}
S_r(V,\d)=\omega_1^r+\cdots+\omega_i^r-\omega_{i+1}^r-\cdots -\omega_\kappa^r,
\end{equation}
où $\omega_1,\ldots,\omega_i$ sont les zéros du numérateur, $\omega_{i+1},\ldots , \omega_\kappa$ ceux du dénominateur de $L$. 
 Deligne a montré que
pour chaque $j$, $|\omega_j|=p^{m_j/2}$ où $m_j\in\N$.

\medskip

Ainsi, pour démontrer le Théorème \ref{expo}, puisque
$$S(\b,\cc,\d)=S_1^*(V, \d)\;,$$
il suffit de prouver que $\kappa\leq 44^{4k-1}$ et $m_j\le 2$ pour tout $1\le j\le \kappa$, ce que nous allons faire maintenant.

\subsection{Majoration du nombre  $\kappa$ de la formule \eqref{SrVd}}

Notre point de départ est le résultat suivant

\begin{theo}\label{kappaH} {\bf (Hooley\cite{H82}  Theorem 4 p. 112)}
Soit $f\in\F_p [X_1,\ldots ,X_N]$ de degré $d$.
Pour $n\in\N$, on définit $$S_n(f)=\sum_{x\in\F_{p^n}^N}\e_p (\sigma_n (f(x))\ {\rm et}\
S_n^*(f)=\sum_{x\in(\F_{p^n}^*)^N}\e_p (\sigma_n (f(x)).$$
Les sommes d'exponentielles $S_n(f)$ et $S_n^*(f)$
admettent une écriture sous la forme \eqref{SrVd} et le nombre des pôles
correspondant $\kappa$ vérifie 
$$\kappa \le (11d+11)^{N+1}.$$
\end{theo}

En utilisant l'orthogonalité des caractères (comme l'a fait Hooley (\cite{H82} p. 104)  pour détecter les conditions définissant $V_n$, on remarque que
\begin{equation}
S_n(V,\d )=\frac{1}{p^{2(k-1)n}}S_n (\phi),
\end{equation}
avec pour $u =(u_1,\ldots ,u_k)$, $v=(v_1,\ldots ,v_k)$, $g=(g_2,\ldots ,g_k)$, $h=(h_2,\ldots, h_k)$ :
$$
\phi ( u,v, g, h)=d_1u_1v_1+
\sum_{i=2}^k (d_iu_iv_i
+g_i(u_iu_1(b_1-b_i)
+u_i-u_1)+h_i(v_iv_1(c_1-c_i)+v_i-v_1)).
$$
D'après le Théorème \ref{kappaH} (avec $N=4k-2$ et $d=3$),
le $\kappa$ correspondant à  $S_n(V,\d)$ est inférieur à  $44^{4k-1}$.

\subsection{Majorations des puissances $m_j$ de \eqref{SrVd}}
Nous allons montrer que $m_j\le 2$ pour tout $1\le j\le \kappa$.

On procède comme Hooley \cite{H82} ou comme \cite{BiBo85} avec un argument de valeur moyenne.

Pour $\lambda\in\F_{p^n}$, on considère la variété

$$W_\lambda=\{ (x,y)\in V_n ,\ \sum_{i=1}^kd_ix_iy_i=\lambda\}.$$
Il sera parfois utile de noter $W_\lambda (\F _{p^n})$ cette variété.
On a vu que la définition de $V_n$ implique que  si l'un des $x_i$ est nul  alors $x_1=\ldots =x_k=0$.
Ainsi, pour $\lambda\not =0$, 
$$W_\lambda=\{ (x,y)\in V_n ,\ \sum_{i=1}^kd_ix_iy_i=\lambda, \prod_{i=1}^k x_i y_i\not = 0\}.$$
Remarquons ({\it via} la correspondance $x_i=\ov{x-b_i}$ et $y_i=\ov{y-c_i}$) l'égalité pour $\lambda\not =0$ :

$$W_\lambda=\{ (x,y)\in \F_{p^n}^2 : \sum_{i=1}^kd_i\ov{(x-b_i)(y-c_i)}=\lambda, \prod_{i=1}^k(x-b_i)(y-c_i)\not =0\}.$$
(Pour $\lambda =0$ il faut rajouter le point $(0,0)$.)

Soit $N_n(\lambda )$ le nombre de points de $W_\lambda$.

On reprend maintenant pas à  pas les arguments de Hooley \cite{H82}  ou de Birch et Bombieri \cite{BiBo85}
qui consistent à  évaluer de deux manières différentes la quantité
$$S:=\sum_{c\in \F_{p}\setminus\{ 0\}} | S_c|^2,$$
avec 
$$S_c=
 \sum_{(x,y)\in V_n}\e_p\Big (\sigma_n (c\sum_{i=1}^kd_ix_iy_i)\Big ).$$
Tout d'abord, en regroupant les $x,y$ tels que $\sum_{i=1}d_ix_iy_i=\lambda$, on remarque que les sommes
$S_c$ vallent : 
$$S_c =\sum_{\lambda\in\F_{p^n}}N_n(\lambda)
\e_p(\sigma_n(c\lambda)).$$
On obtient en développant les carrés de la somme $S$
\begin{displaymath}
\begin{array}{ll}
S&=\sum_{c\in \F_{p}\setminus\{ 0\}} \sum_{\lambda,\lambda '\in\F_{p^n}}
N_n(\lambda )N_n(\lambda ')\e_p (\sigma_n(c\lambda)-\sigma_n(c\lambda '))\\
&=p^n\sum_{\lambda\in\F_{p^n}}N_n(\lambda )^2-(\sum_{\lambda\in\F_{p^n}}N_n(\lambda ))^2,
\end{array}
\end{displaymath}
soit
$$S=p^n\sum_{\lambda\in\F_{p^n}}(N_n(\lambda )-M)^2,$$
où $M$ est la valeur moyenne :
$$M=\frac{1}{p^n}\sum_{\lambda\in\F_{p^n}}N_n(\lambda ) =p^n.$$
Posons $$g_\lambda (X,Y)=\sum_{i=1}^k d_i\prod_{j\not = i} (X-b_j)(Y-c_j)-\lambda\prod_{i=1}^k(X-b_i)(Y-c_i),$$
de sorte que pour $\lambda\not =0$, $$W_\lambda (\F _{p^n})=\{ (x,y)\in\F_{p^n}^2 : g_\lambda (x,y)=0, 
\prod_{i=1}^k(x-b_i)(y-c_i)\not =0\}.$$
Pour chaque $x\in\F _p^n$, $y\mapsto g_\lambda (x,y)$, est la fonction polynomiale d'un polynôme de
$\F_{p^n}[Y]$ de degré $k-1$ ou $k$ suivant que $x$ soit l'un des $b_i$ ou non.
Il a donc dans $\F_{p^n}$ au plus $k$ racines.
On en déduit la première majoration triviale :
\begin{equation}
\label{Nlambda}
N_n(\lambda )= |W_\lambda (\F_{p^n})|\le kp^n.
\end{equation}
Nous obtenons maintenant une expression plus précise pour presque tout $\lambda$ :

\begin{proposition}
\label{lesvlambda}
Il existe deux constantes $\nu,K$ telles pour presque tout $\lambda\in\F_{p^n}$ avec au plus $K$ exceptions, on a :

\begin{equation}
\label{Nn}|N_n(\lambda )-p^n|\le \nu\sqrt{p^n}.
\end{equation}
De plus, les valeurs $\nu =2k^2$ et $K = 9(k-1)^2+4(k-1)+2$ sont admissibles. 
\end{proposition}
La preuve de cette proposition n'est pas immédiate. Avant de l'exposer nous proposons de montrer que cette proposition
est suffisante pour démontrer le Théorème \ref{expo}.
Admettons donc provisoirement la Proposition \ref{lesvlambda}.
Nous reprenons les idées de Hooley \cite{H82} pp. 115-116. 
Comme nous voulons contrôler la dépendance en $k$, nous les reproduisons ici dans notre contexte.
 Il s'agit notamment de ne pas rater  
un ``assez grand'' qui dépendrait de $k$.
Soit $\cK$ l'ensemble des $\lambda \in \F_{p^n}$ ne vérifiant pas la Proposition \ref{lesvlambda}.
 On a alors

\begin{equation}
\label{calS1}S=p^n\sum_{\lambda \in\cK} (N_n(\lambda )-p^n)^2+p^n\sum_{\lambda\not\in\cK} (N_n(\lambda)-p^n)^2
\le (Kk^2+\nu ^2) p^{3n}.
\end{equation}

Supposons que  l'un des $\omega_i$ de la formule \eqref{SrVd} soit de module $p^{m/2}$ pour un $m\ge 3$.
 Hooley a montré à partir   des travaux de Deligne que les sommes 
$S_c$ pour $c\in\F_p^*$ sont également de la forme 
$$S_c=e_1\omega_{1,c}^n+\cdots +e_L\omega_{L,c}^n,$$
où $e_1,\ldots, e_L$ sont des entiers indépendants de $c$ et $n$ vérifiant $|e_1|+\cdots +|e_L|=\kappa$,
et  les modules $|\omega_{j,c}|$ sont des puissances entières de $\sqrt{p}$ ;  ces puissances étant indépendantes de $c$.

En particulier le nombre des $\omega _{i,c}$ de module supérieur à $p^{3/2}$ est le même pour chaque $c$.
Plus précisément si $\ell$ désigne le nombre d'indices $i$ tels que $\omega_{i,c}$ soit de module supérieur 
à $p^{3/2}$, on a  pour tout $c\in\F_p^*$ (quitte à changer l'ordre des $\omega_{i,c}$) :
$$S_c =e_1\omega_{1,c}^n+\cdots +e_\ell\omega_{\ell ,c}^n+E_c,$$
où $E_c$ est un terme d'erreur  de module inférieur à $\kappa p^{n}.$
Soit $H$ tel que $p^H=\max_{1\le i\le L} |\omega _{i,c}|^2$.  
Posons $z_{i,c}=\frac{\omega_{i,c}}{p^{H/2}}$ pour $1\le i\le \ell$. Les $z_{i,c}$ sont ainsi des nombres complexes deux à deux distincts de modules inférieurs où égaux à $1$.
On en déduit la minoration  pour tout $c\in\F_p^*$ :

$$|S_c|\ge p^{3n/2}|e_1z_{1,c}^n+\cdots +e_\ell z_{\ell ,c}^n|-\kappa p^n.$$
Cela donne pour $S$ :

$$S\ge \sum_{c\in\F_p^*}|S_c|^2\ge \sum_{c\in\F_p^*}p^{3n}|e_1z_{1,c}^n+\cdots +e_\ell z_{\ell ,c}^n|^2-2\kappa ^2(p-1)p^{5n/2}.$$

L'idée suivante de Hooley est d'utiliser l'égalité :
$$\lim_{u\rightarrow +\infty}\frac{1}{u}\sum_{n\le u}|e_1z_{1,c}^n+\cdots +e_\ell z^n_{\ell ,c}|^2=|e_1^2+\cdots +e_\ell ^2|
\ge 1,$$
qui se vérifie en développant le carré et en profitant du fait que les $z_{i,c}$ sont des nombres complexes deux à deux distincts de module inférieur à $1$.
Grâce à cette égalité, on observe que pour $\varepsilon >0$ donné, il existe une infinité d'entiers $n$ tels que 
\begin{equation}
\label{grosn}
|e_1z_{1,c}^n+\cdots +e_\ell z_{\ell ,c}^n|^2\ge 1-\varepsilon.
\end{equation}

On obtient alors pour les entiers $n$ vérifiant \eqref{grosn} :

$$ S\ge (p-1)p^{3n} \big ( 1-\varepsilon-\frac{2\kappa ^2}{p^{n/2}}\big )\ge (p-1)p^{3n}(1-2\varepsilon),$$

pour $n$ assez grand.
Cette minoration est incompatible avec \eqref{calS1} lorsque $p>Kk^2+\nu ^2 +1$ (avec un choix de $\varepsilon$ assez petit).

Lorsque $p\le Kk^2+\nu ^2+1$, le Théorème \ref{expo} reste vrai mais fournit en fait une majoration moins bonne que la majoration triviale ($|S(\b,\cc ,\d)|\le p^2$).
\medskip

Il reste maintenant à  démontrer la proposition \ref{lesvlambda} pour terminer la preuve du Théorème
\ref{expo}.

Le cardinal $N_n(\lambda )$ est proche du cardinal
$$N'_n(\lambda ) :=|\{ (x,y)\in\F_{p^n}^2 : g_\lambda (x,y) =0\}|.$$
Pour obtenir une majoration de la différence $|N'_n(\lambda) -N_n(\lambda )|,$ il suffit  d'étudier la contribution des $x=b_i$ ou $y=c_i$ dans $g_\lambda$.
Le polynôme associé à
$$g_\lambda (b_i,y)=d_i\prod_{\substack{1\leq j\leq k \\ j\not = i}} (b_i-b_j)(y-c_j)$$
 a  $k-1$ racines. De même, pour $1\le i\le k$,
 $|\{ x\in\F_{p^n}: g_\lambda (x,c_i)=0\}|= k-1$.
Ainsi,
\begin{equation}
\label{Nprime}|N'_n(\lambda) -N_n(\lambda )|\le 2k.
\end{equation}
Pour démontrer la Proposition \ref{lesvlambda}, il suffit donc d'étudier $N'_n(\lambda )-p^n$.

Si $g_\lambda$ est absolument irréductible (c'est-à -dire irréductible dans $\ov{\F_p}$), alors d'après un théorème de Lang et Weil \cite{LW54}, $N'_n(\lambda )=p^n+O(\sqrt{p^n})$.
Nous ne savons pas dire que $g_\lambda$ est absolument irréductible pour tout $\lambda$ sauf pour au plus un nombre fini de $\lambda$.

Nous utilisons plutôt les travaux d'Adolphson et Sperber \cite{AS89}.
Pour cela nous devons définir
le polygone de Newton d'un polynôme $g\in K[X_1,\ldots ,X_n]$ où $K$ est un corps fini.
\'Ecrivons $g$ sous la forme $g=\sum_{j\in J}a_jx^j$ où $J\subset\Z_+^n$,
et pour chaque $j=(j_1,\ldots ,j_n)\in J$, $a_j\not =0$ et
$x^j=x_1^{j_1}\cdots x_n^{j_n}$.

Le polygone de Newton de $g$ est alors l'enveloppe convexe dans $\R^n$ de l'ensemble $J\cup\{ (0,\ldots ,0)\}$. On le note $\Delta (g)$.
La dimension de $\Delta (g)$ est celle du plus petit sous-espace vectoriel de $\R^n$ contenant $\Delta (g)$.

\`A chaque face $\sigma$ de $\Delta (g)$, on associe le polynôme
$$g_\sigma=\sum_{j\in\sigma\cap J}a_jx^j.$$

\begin{definition}
\label{degcom}
(i) Le polynôme $g$ est non dégénéré (par rapport à  $\Delta (g)$) si pour toute
face $\sigma$ de $\Delta (g)$ qui ne contient pas l'origine, les polynômes
$\partial g_\sigma/\partial x_1$,\ldots, $\partial g_\sigma/\partial x_n$ n'ont pas de racine commune dans $(\ov{K^*})^n$ où $\ov K$ désigne une clôture algébrique de $K$.

(ii) Le polynôme $g$ est commode (par rapport à  $\Delta(g)$)  si pour tout $i\in \{1,\ldots ,n\}$, il existe un entier $j_i>0$ tel que $g$ contienne
un monôme de la forme $ax_i^{j_i}$.
\end{definition}
Ces deux définitions sont celles de \cite{AS89} page 376 adaptées à  notre situation

Adolphson et Sperber ont obtenu le résultat suivant.

\begin{theo}
\label{AASS} {\bf (Adolphson et Sperber \cite{AS89} Corollary 6.9 p. 400)}
Soient $K$ un corps fini de cardinal $q$ et $g\in K[x_1,\ldots ,x_\ell]$ non dégénéré et commode par rapport à  son polygone de Newton. On suppose aussi que
$$g, x_1\partial g/\partial x_1,\ldots ,x_\ell\partial g/\partial x_\ell$$
n'ont pas de zéro commun.
Soit $V$ la variété définie par $g=0$.
Il existe $\nu (g)$ telle que
$$|V(K)|-q^{\ell-1}|\le \nu (g)\sqrt{q^{\ell-1}}.$$
\end{theo}
Remarque : la constante $\nu (g)$ est effectivement calculable. On trouvera une définition de cette constante à  la page 371 de \cite{AS89}, elle ne dépend que du polygone de Newton de $g$.

Pour terminer la preuve de la Proposition \ref{lesvlambda} il reste  à  vérifier que l'on peut appliquer le Théorème \ref{AASS} au polynôme $g_\lambda$ avec $K=\F_{p^n}$ et $\ell=2$ pour presque tous $\lambda$ avec au plus $K$ exceptions
et déterminer les $\nu (g_\lambda )$ correspondants.
Rappelons la forme de $g_\lambda$ :

$$g_\lambda (X,Y)=\sum_{i=1}^kd_i\prod_{j\not =i}{(X-b_j)(Y-c_j)}-
\lambda\prod_{i=1}^k (X-b_i)(Y-c_i).$$

Lorsque $\lambda \not = 0$, le coefficient en $X^k$ vaut
$-\lambda (-1)^kc_1\ldots c_k\not =0$ puisque aucun $c_i$ n'est nul.

On vérifie de la même manière pour $\lambda\not =0$ que le coefficient
en $Y^k$ est non nul.

Cela prouve que $g_\lambda$ est commode si $\lambda\not = 0$.

Le polygone de Newton est $\Delta (g_\lambda )=[0,k]\times [0,k]$ si $\lambda
\not =0$.

En reprenant la définition de $\nu$ donnée page 371 de \cite{AS89} (avec $\nu =\nu_A$ où $A=\{ 1,2\}$)
on obtient sans peine que
\begin{equation}
\label{nug}
\nu (g_\lambda )=2k^2-2k.
\end{equation}

Les faces de $\Delta (g_\lambda)$ ne contenant pas $(0,0)$ sont les cotés $\sigma_1=\{k\}\times[0,k]$ et $\sigma_2=[0,k]\times\{ k\}$.

Les polynômes associés sont $g_{\sigma_1}(X,Y)=-\lambda X^k\prod_{i=1}^k
(Y-c_i)$ et $g_{\sigma_2}(X,Y)=-\lambda Y^k\prod_{i=1}^k (X-b_i).$

Vu que les $c_i$ sont deux à  deux distincts, on vérifie facilement
que $\partial g_{\sigma _1}/\partial x$ et $\partial g_{\sigma_1}/\partial y$ n'ont pas de racine commune dans $(\ov{\F_{p^n}^*})^2$.
Cette propriété est également vérifiée par $g_{\sigma _2}$ car les $b_i$ sont deux à  deux distincts.
Donc $g_\lambda$ est non dégénéré lorsque $\lambda\not =0$.

La dernière condition du Théorème \ref{AASS} est plus difficile à  vérifier.

Commençons par traiter le cas où $x=b_i$ pour un $i\in\{1,\ldots ,k\}$.
Alors
$$g_\lambda (b_i,Y)=d_i\prod_{j\not =i}(b_i-b_j)(Y-c_j)$$
s'annule pour $Y=c_j$ avec $j\not =i$.
Mais
$$\frac{\partial g_\lambda}{\partial x}(b_i,c_j)=d_j(c_j-c_i)\prod_{\ell\not =i,j}(b_i-b_\ell)(c_j-c_\ell)\not =0.$$
Donc pour tout $\lambda$, il n'existe pas de point singulier de la forme $(b_i,y)$.
De même, il n'existe pas de point singulier de la forme $(x,c_i)$.

Pour $x\not\in\{b_1,\ldots ,b_k\}$, $y\not\in\{ c_1,\ldots c_k\}$, on vérifie que
$g_\lambda (x,y)=0$ si et seulement si
 \begin{equation}\label{lambdadeg}\lambda=\sum_{i=1}^kd_i\ov{x-b_i}\ov{y-c_i}.
\end{equation}
Pour la dérivée partielle en $y$, cela entraîne :
\begin{equation}\label{party}
\begin{split}
\frac{\partial g_\lambda}{\partial y} (x,y)&=\sum_{i=1}^kd_i\prod_{j\not =i}(x-b_j)
\sum_{j\not =i}\prod_{\ell\not = i,j}(y-c_\ell)-\lambda\prod_{i=1}^k(x-b_i)\sum_{i=1}^k\prod_{j\not =i}(y-c_j)\\
&=-\prod_{i=1}^k(x-b_i)(y-c_i)\sum_{i=1}^kd_i\ov{(x-b_i)(y-c_i)^2}.\\
\end{split}
\end{equation}
De même,
$$\frac{\partial g_\lambda}{\partial x}(x,y)=-\prod_{i=1}^k(x-b_i)(y-c_i)\sum_{i=1}^kd_i\ov{(x-b_i)^2(y-c_i)}.$$

Tout d'abord on considère le cas où $xy=0$. Pour $(x,y) =(0,0)$ le système devient
\begin{equation}
\label{00} g_\lambda (0,0)=0\Leftrightarrow \sum_{i=1}^kd_i\prod_{j\not =i}
b_jc_j-\lambda\prod_{i=1}^kb_ic_i=0,
\end{equation}
ce qui n'arrive que pour  une seule valeur de $\lambda$ puisque les $b_i$ et les $c_i$ ne sont pas nuls.

Lorsque $x=0$ et $y\not =0$, on doit résoudre le système
\begin{equation}\label{x0}
g_\lambda (0,y)=\frac{\partial g_\lambda}{\partial y} (0,y)=0.
\end{equation}

La deuxième équation devient (rappelons que $y\not\in\{ c_1,\ldots ,c_k\}$)

$$\sum_{i=1}^kd_i\ov{b_i(y-c_i)^2}=0$$
ou encore en multipliant par $\prod_{i=1}^k(y-c_i)^2$ :
$$\sum_{i=1}^kd_i\ov{b_i}\prod_{j\not =i}(y-c_j)^2=0.$$
Comme les $c_j$ sont deux à  deux distincts et les $b_i$ ne sont pas nuls
ce polynôme ne s'annule pas en $y=c_1$. Il  n'est donc pas identiquement nul, et admet au plus $2(k-1)$
racines $y$.
On en déduit en utilisant \eqref{lambdadeg} qu'il n'y a au plus que $2(k-1)$ valeurs de  $\lambda$ telles que le
système \eqref{x0} ait des solutions.

On vérifie de la même façon que le système
\begin{equation}\label{y0}
g_\lambda (x,0)=\frac{\partial g_\lambda}{\partial x} (x,0)=0
\end{equation}
admet des solutions pour au plus $2(k-1)$ valeurs de $\lambda$.

Il reste maintenant à  étudier le système
\begin{equation}\label{xy}
g_\lambda (x,y)= \frac{\partial g_\lambda}{\partial y} (x,y)= \frac{\partial g_\lambda}{\partial y} (x,y)=0,
\end{equation}
avec   $x\not\in\{ b_1,\ldots ,b_k\}$, $y\not\in\{ c_1,\ldots ,c_k\}$.

En reprenant le calcul précédent on remarque que ce système est équivalent à  :
\begin{equation}\label{xypanul}
\left\{
\begin{split}
\sum_{i=1}^kd_i\ov{(x-b_i)(y-c_i)} &=\lambda \\
\sum_{i=1}^kd_i\prod_{j\not = i}(x-b_j)^2(y-c_j) &=0\\
 \sum_{i=1}^kd_i\prod_{j\not = i}(x-b_j)(y-c_j)^2 &=0
\end{split}\right.
\end{equation}
Soit $P_1(x,y)=\sum_{i=1}^kd_i\prod_{j\not = i}(x-b_j)^2(y-c_j) $
et $P_2(x,y)=\sum_{i=1}^kd_i\prod_{j\not = i}(x-b_j)(y-c_j)^2.$

Pour montrer que les deux dernières lignes du système n'ont qu'un nombre fini de solutions $(x,y)$ il suffit de vérifier que les polynômes $P_1$ et $P_2$ sont premiers entre eux.
Notons $T=(P_1,P_2)$ puis $P_1=TR_1$ et $P_2=TR_2$.
Si $T$ n'est pas constant alors quitte à  échanger les rôles de $x$ et de $y$, on peut supposer que le degré  partiel en $y$ pour $T$ est supérieur à  $1$.

Le degré en $y$ de $P_1$ est inférieur à  $k-1$ c'est donc aussi le cas
pour $T$ et $R_1$. De plus,  le degré partiel en $y$ de $R_1$ est au plus $k-2$.
Comme les $c_j$ sont deux à  deux distincts, les polynômes
$\prod_{j\not =i}(y-c_j)$ forment une base de l'espace vectoriel des polynômes de $\F_{p^n}[Y]$ de degré au plus $k-1$.

On en déduit que $T$ et $R_1$ s'écrivent sous la forme :
\begin{equation}\label{R1T}
T(X,Y)=\sum_{i=1}^k\alpha _i(X)\prod_{j\not =i}(Y-c_j),\quad
R_1(X,Y)=\sum_{i=1}^k\beta _i(X)\prod_{j\not =i}(Y-c_j),
\end{equation}
où $\alpha _i (X),\beta _i(X)\in\F_{p^n}[X]$ pour $1\le i\le k$.
Comme $T$ divise $P_2$, pour tout $1\le i\le k$ le polynôme
$T(X,c_i)=\alpha _i (X)\prod_{j\not =i}(c_i-c_j)$ divise $P_2(X,c_i)=d_i\prod_{j\not =i}(X-b_j)(c_i-c_j)^2.$
On en déduit que $\alpha_i (X)$ est de la forme
$\alpha _i (X)=s_i\prod_{j\in L_i}(X-b_j)$ avec $s_i\in\F_{p^n}$
et $L_i\subset \{ 1,\ldots ,k\}\setminus\{ i\}$.

On a alors
$$P_1(X,c_i)=d_i\prod_{j\not =i}(X-b_j)^2(c_i-c_j)=s_i\prod_{j\in L_i}(X-b_j)
\beta_i (X)\prod_{j\not =i}(c_i-c_j)^2.$$
Ainsi le polynôme $\beta_i$ est de la forme
$$\beta _i (X)=t_i\prod_{j\not =i}(X-b_j)\prod_{j\not\in L_i\cup\{i\}}
(X-b_j).$$

Le coefficient du terme en $Y^{k-1}$ dans l'écriture de $R_1$ dans \eqref{R1T} est $\sum_{i=1}^k\beta_i (X)$. Ce coefficient doit être nul car
$R_1$ est de degré au plus $k-2$ en Y.
Or, pour $1\le i\le k$,
$$0=\sum_{m=1}^k\beta _m(b_i)=\beta _i(b_i)=t_i\prod_{j\not =i}(b_i-b_j)\prod_{j\not\in L_i\cup\{ i\}}(b_i-b_j).$$
Cela entraîne que $t_i=0$ puisque les $b_j$ sont deux à  deux distincts.
Ainsi, chaque $\beta_i$ et par suite $R_1$ et $P_1$ sont identiquement nuls.
Cela n'est pas possible.
Par conséquent $P_1$ et $P_2$ sont premiers entre eux.

Les polynômes homogènes associés $P_1(X:Y:Z)=Z^{3(k-1)}P_1(X/Z,Y/Z)$, $P_2(X:Y:Z)=Z^{3(k-1)}P_2(X/Z,Y/Z)$
sont alors premiers entre eux. En effet leur pgcd $h$ est un polynôme homogène  vérifiant $h(X:Y:1)=1$ et est donc de la forme 
$h(X:Y:Z)=aZ^t$ ce qui n'est possible que pour $t=0$. 
On en déduit par un résultat classique  sur les intersections de courbes planes (Théorème de Bézout,  cf \cite{Wald} ou \cite{Walk})
que 
$$|\{  (x,y)\in\F_{p^n}^2 : P_1(X,Y)=0=P_2(X,Y)\}|\le \deg P_1 \times \deg P_2 \le 9(k-1)^2.$$

Le Théorème \ref{AASS} s'applique donc pour presque tous les $\lambda$ avec au plus 
$K=9(k-1)^2+4(k-1)+2$ exceptions. On rappelle que le terme $4(k-1)+1$ correspond aux $\lambda$ 
pour lesquels il existe un couple $(x,y)$ tels que $xy=0$ et
solutions d'un des systèmes \eqref{00}, \eqref{x0} ou \eqref{y0} et le $+1$ supplémentaire tient compte du cas $\lambda =0$.
En tenant compte de \eqref{Nprime} on obtient la Proposition
\ref{lesvlambda} et cela termine la preuve du  Théorème \ref{expo}.

\section{La complexité $K_3$ dans le cas où $S$ est l'ensemble des inverses d'une suite d'entiers consécutifs}
Dans ce paragraphe nous démontrons le théorème \ref{thm:K3}
relatif aux polynômes à  racines simples dans $\F _p$.
Pour $\cA\subset\F _p$, on reprend la notation $f_\cA (X)=\prod_{a\in\cA}(X-a).$

On remarque que  $|\cP_3 (d,p)|={p\choose d}$. La condition $p>2/(1-\beta )$ implique que $p\ge 3$.

Soient $\cB$ et $\cC$ deux sous-ensembles disjoints de $\F _p$ :
$$\cB=\{ b_1,\ldots , b_\ell\},\quad\cC=\{ b_{\ell+1},\ldots ,b_k\}.$$

Nous devons montrer qu'il existe $f\in\cP _3 (d,p)$ tel que
$r_p(\ov{f(b_i)})\in \{ r,\ldots, r+s-1\}$ pour $1\le i\le\ell$
et
$r_p (\ov{f(b_i)})\not\in \{ r,\ldots, r+s-1\}$ pour $\ell <i\le k$.
Ici et dans la suite nous considérons des polynômes $f=f_\cA$ tels que $\cA\cap (\cB\cup\cC)=\emptyset$.

Pour $\beta\in ]0,1[$ on définit la fonction $h_\beta (x)$ pour $x\in\Z$ par
\begin{displaymath}
h_\beta (x)=\left\{ \begin{array}{ll}  1 & \text{si}\ 0\le r_p (x)<\beta p\\
0 & \text{sinon}.
\end{array}\right.
\end{displaymath}
En adaptant la preuve du lemme 2.2 de \cite{DMS}, on montre que
\begin{equation}\label{hbeta}h_\beta (x)=\sum_{|h|<p/2} \alpha_h \e _p (hx),
\end{equation}
avec
$$\alpha_0=\frac{\rf{\beta p}}{p}\ {\rm et} \ \alpha_h =\frac{1-\e_p(-h\rf{\beta p})}{p(1-
\e_p (-h))} \ {\rm pour}\ h\not= 0.$$
On a ensuite remarqué dans \cite{DMS} que $|\alpha_h|\le 1/(2h)$ pour $h\not =0$.

Il suffit de montrer que la quantité
$$T:=\sum_{\substack{f\in\cP_3(d,p)\\ \prod_{i=1}^kf(b_i)\not =0}}\prod_{i=1}^\ell h_\beta (\ov{f(b_i)}-r)\prod_{i=\ell +1}^k
(1-h_\beta (\ov{f(b_i)}-r))$$
est strictement positive.
On développe les produits ci-dessus en utilisant \eqref{hbeta} et en isolant le terme principal :
\begin{displaymath}
\begin{split}
T&= {p-k\choose d}\frac{\rf{\beta p}^\ell(p-\rf{\beta p})^{k-\ell}}{p^k}\\
&+\sum_{\substack{0\le t\le\ell\\ 0\le t'\le k-\ell\\ (t,t')\not = (0,0)}}
\frac{\rf{\beta p}^{\ell -t}(p-\rf{\beta p})^{k-\ell -t'}}{p^{k-t-t'}}
(-1)^{t'}\\
&\times
\sum_{\substack{1\le i_1<\ldots <i_t\le\ell\\ \ell< i_{t+1}<\ldots <i_{t+t'}
\le k}}
\sum_{0<|h_1|,\ldots ,|h_{t+t'}|<p/2}\Big (\prod_{i=1}^{t+t'}\alpha_{h_i}\Big )
\e_p \Big ( -r\sum_{m=1}^{t+t'}h_m\Big )
U(h_1,\ldots ,h_{t+t'}),\\
\end{split}
\end{displaymath}
avec

\begin{displaymath}
\begin{split}
U(h_1,\ldots ,h_{t+t'}) &=\sum_{\substack{f\in\cP _3(d,p)\\  \prod_{i=1}^kf(b_i)\not =0}}
\e_p\Big (
\sum_{m=1}^{t+t'}h_m\ov{f(b_{i_m})}\Big )\\
&=\sum_{ \{ a_1,\ldots ,a_d\}\subset\F_p\setminus (\cB\cup\cC)}
\e _p \Big ( \sum_{m=1}^{t+t'}h_m\prod_{j=1}^d \ov{ b_{i_m}-a_j}
\Big ).\\
\end{split}
\end{displaymath}

On fixe les $d-2$ premières variables et on somme sur les deux dernières :

\begin{displaymath}
\begin{split}
|U(h_1,\ldots ,h_{t+t'}) | &\le \sum_{ \{ a_1,\ldots ,a_{d-2}\}\subset\F_p\setminus (\cB\cup\cC)}
\left | \sum_{\substack{\{ a_{d-1}, a_d\}\subset\F_p\setminus \cB\cup\cC\\
a_i\not = a_j\ si \ i\not =j}}
\e _p \Big ( \sum_{m=1}^{t+t'}h_m\prod_{j=1}^d \ov{ b_{i_m}-a_j}
\Big )\right |\\
&\le \sum_{ \{ a_1,\ldots ,a_{d-2}\}\subset\F_p\setminus (\cB\cup\cC)}
\left | \sum_{a_{d-1},a_d\in \F_p\setminus (\cB\cup\cC )}
\e _p \Big ( \sum_{m=1}^{t+t'}h_m\prod_{j=1}^d \ov{ b_{i_m}-a_j}
\Big )\right |\\
&+ 2(k+d-1)p^{d-1}.
\end{split}
\end{displaymath}
Ce terme $2(k+d-1)p^{d-1}$ est une majoration de la contribution des $(a_1,\ldots ,a_d)$
tels que $a_d=a_{d-1}$, ou $a_{d-1},a_d\in\cB\cup\cC\cup\{a_1,\ldots ,a_{d-2}\}$.
D'après le Théorème \ref{expo} la somme intérieure sur $a_{d-1}, a_d$ est majorée par
$44^{4k}p$.
On en déduit la majoration :
$$|U(h_1,\ldots ,h_{t+t'}) |\le (44^{4k}+2k+2d-2)p^{d-1}.$$
Cela donne pour $T$ (la constante implicite dans la formule suivante est de valeur absolue inférieure à  1):
\begin{displaymath}
\begin{split}
T&= {p-k\choose d}\frac{\rf{\beta p}^\ell(p-\rf{\beta p})^{k-\ell}}{p^k}\\
&+O\Big (\sum_{\substack{0\le t\le\ell\\ 0\le t'\le k-\ell\\ (t,t')\not = (0,0)}}
\frac{\rf{\beta p}^{\ell -t}(p-\rf{\beta p})^{k-\ell -t'}}{p^{k-t-t'}}(-1)^{t+t'}\\
&\times
\sum_{\substack{1\le i_1<\ldots <i_t\le\ell\\ \ell< i_{t+1}<\ldots <i_{t+t'}
\le k}}
\sum_{0<|h_1|,\ldots ,|h_{t+t'}|<p/2}\Big (\prod_{i=1}^{t+t'}|\alpha_{h_i}|\Big )
 (44^{4k}+2k+2d-2)p^{d-1}\Big ).
\end{split}
\end{displaymath}
Comme $|\alpha_h|\le (2h)^{-1}$ lorsque $h\not =0$, les sommes sur les $h_i$
sont inférieures à  $\log (3p)$.
On obtient donc :
\begin{displaymath}
\begin{split}
T &= {p-k\choose d}\frac{\rf{\beta p}^\ell(p-\rf{\beta p})^{k-\ell}}{p^k}\\
&+O  \Big ( \Big ( \frac{\rf{\beta p}}{p}+\log (3p)\Big )^\ell
\Big ( \frac{p-\rf{\beta p}}{p}+\log (3p)\Big )^{k-\ell}
(44^{4k}+2k+2d-2)p^{d-1}.\\
&= {p-k\choose d}\frac{\rf{\beta p}^\ell(p-\rf{\beta p})^{k-\ell}}{p^k}
+O\Big ( \log (9p)^k(44^{4k}+2k+2d-2)p^{d-1}\Big ),
\end{split}
\end{displaymath}
où la constante implicite du $O$ est inférieure à  $1$ en valeur absolue.
Lorsque $p$ vérifie \eqref{condK3p} $T>0$ et ainsi $K_3(S,d)\ge k$.
Cela termine la démonstration du Théorème \ref{thm:K3}.

\section{La complexité $K_3$ dans le cas général}
Dans ce paragraphe on étudie
$K_3 (S,d)$ où $S$ est maintenant un sous-ensemble de $\F_p$ vérifiant les hypothèses du Théorème \ref{k3general}.

Soit $k\le d/2$. Soient $\cB,\cC$ deux sous-ensembles de $\F_p$ disjoints tels
que $|\cB | =\ell, |\cC |=k-\ell$ pour un certain
$0\le\ell \le k$.
Nous devons montrer qu'il existe un sous-ensemble $\cA\subset\F_p$
de cardinal $d$ tel que $\cB\subset \cR (f_\cA, S)$ et $\cC\cap \cR (f_\cA, S)=\emptyset$.

Notons $\cB=\{b_1,\ldots ,b_\ell\}$, $\cC=\{ b_{\ell +1},\ldots,b_k\}$.

\medskip

$\bullet$ Si $S=\{ 0\}$, alors on peut prendre
$f_\cA (X)=\prod_{i=1}^\ell (X-b_i)\prod_{i=\ell +1}^d(X-b'_i)$
avec $\{ b'_{\ell +1},\ldots ,b'_d\}\cap (\cB\cup\cC)=\emptyset$,
ce qui est possible lorsque $k+d<p$.
De même pour $S^c=\{ 0\}$,
$f_\cA (X)=\prod_{i=1}^{d-\ell -1} (X-b'_i)\prod_{i=\ell +1}^k(X-b_i)$
avec $\{ b'_1,\ldots ,b'_{d-\ell-1}\}\cap (\cB\cup\cC)=\emptyset$ convient.

\medskip

$\bullet$ Nous supposons maintenant que $S\not =\{0\}, \F _p^*$.
Il existe $s\in S\setminus\{0\}$
et $r\in S^c\setminus\{ 0\}$.

Pour terminer la preuve du théorème \ref{k3general}, il suffit de trouver $\cA\subset \F_p$ de cardinal $d$ tel que

\begin{equation}\label{rs}
f_\cA (x)\in\left\{ \begin{array}{ll}  S\setminus\{ 0\} & \text{si} \ x\in\cB\\
S^c\setminus\{ 0\} & \text{si}\ x\in\cC.
\end{array}\right.
\end{equation}

Soit $T_2(S,d)$ le nombre de sous-ensembles $\cA$ vérifiant \eqref{rs}.

On a alors d'après l'orthogonalité des caractères sur $\F _p$ :

\begin{displaymath}
\begin{split}
\varphi (p)^kT_2 (S,d)& =\sum_{\{a_1,\ldots ,a_d\}\subset \F _p}
\prod_{j=1}^\ell \Big (\sum_{\chi }\sum_{r\in S\setminus\{ 0\}}\chi\big (\prod_{i=1}^d(b_j-a_i)\big )\ov{\chi} (r)\Big )\\
&\times\prod_{j=\ell +1}^k \Big (\sum_{\chi }\sum_{s\in S^c\setminus\{ 0\}}\chi\big (\prod_{i=1}^d(b_j-a_i)\big )\ov{\chi} (s)\Big ),
\end{split}
\end{displaymath}

où dans les différentes sommes, $\chi$ parcourt l'ensemble des caractères de $\F_p^*$, et $\ov{\chi}$ est le caractère conjugué de $\chi$.
A priori, la somme devrait être restreinte aux sous-ensembles
$\cA =\{a_1,\ldots,a_d\}\subset\F_p$ tels que $\cA\cap(\cB\cup\cC)=\emptyset$
mais la contribution des $a_i=b_j$ étant nulle, nous pouvons
oublier cette condition.
Soit $\chi _0$ le caractère principal de $\F _p$.
En développant la ligne ci-dessus et en isolant le terme où tous les
$\chi$ valent $\chi _0$, on obtient :
\begin{displaymath}
\label{chiO}
\begin{split}
\varphi (p)^kT_2 (S,d)& ={p-k\choose d}|S\setminus\{ 0\} |^\ell|S^c\setminus\{ 0\}|^{k-\ell}\\
&+O\Big (\sum_{\substack {0\le t\le \ell\\
0\le t'\le k-\ell\\(t,t')\not =(0,0)}}|S\setminus\{ 0\}|^{\ell -t}|S^c\setminus\{ 0\} |^{k-\ell -t'}\\
&\times\!\!\!\!\!\!\!\!\!\!\!\!\!\sum_{ \substack{1\le i_1<\cdots <i_t\le\ell\\
\ell<i_{t+1}<\cdots<i_{t+t'}\le
k}}\sum_{\chi_1\not =\chi_0}\cdots\sum_{\chi_{t+t'}\not =\chi_0}|Z(\bar \chi_{1},\ldots,\bar \chi_{t+t'})|
|V(\chi_{1},\ldots,\chi_{t+t'})|\Big),
\end{split}
\end{displaymath}
avec
$$V(\chi_1,\ldots,\chi_{t+t'})=\sum_{\cA=\{a_1,\ldots ,a_d\}\subset \F _p}
\prod_{j=1}^{t+t'}\chi_j(f_\cA (b_{i_j})),$$
$$Z(\chi_{1},\ldots,\chi_{t+t'})=\prod_{j=1}^t\Big ( \sum_{r\in S}\chi _j (r)\Big )\prod_{j=t+1}^{t+t'}\Big (\sum_{s\not\in S}\chi _j (s)\Big ),$$
et la constante implicite  est de module inférieur à  $1$.
Pour majorer les sommes $V(\chi_{1},\ldots,\chi_{t+t'})$ nous utilisons
le lemme suivant.

\begin{lemme}
\label{somcar}
 Soient $p$ un nombre premier, $\chi$ un caractère non principal d'ordre $d$ (avec $d| (p-1))$,
$f(X)\in\F_p[X]$. Notons $m$ le nombre de racines distinctes de $f(X)$ dans $\ov{\F_p}$.
Si $f(X)$ n'est pas une puissance $d$ ième alors
$$\left | \sum_{x\in\F_p}\chi (f(x))\right | \le (m-1)\sqrt{p}.$$
\end{lemme}

Il s'agit du théorème 2C' p.43  de \cite{S76} dans le cas où $q=p$.

Avant d'utiliser ce lemme nous devons faire un travail préparatoire analogue à  \cite{DS07b}
 dont les idées de base se trouvent dans \cite{S01}.
La difficulté ici est que nous travaillons avec plusieurs variables
que nous devons rendre
{\it indépendantes}.

Comme $\F_p ^*$ est cyclique, chaque caractère $\chi_m$ peut s'écrire sous la forme
$\chi^{\alpha_m}$ où $\chi$ est un caractère d'ordre $p-1$.
Ainsi pour des caractères $\chi_1,\ldots , \chi_{t+t'}$ différents de $\chi_0$, il existe des entiers  compris entre $1$ et $p-2$,
$\alpha_1,\ldots ,\alpha_{t+t'}$ tels que
$$V(\chi_1,\ldots,\chi_{t+t'})=V(\chi^{\alpha_1},\ldots ,\chi^{\alpha_{t+t'}})=W(\alpha _1,\ldots ,\alpha _{t+t'}),$$
avec
$$W(\alpha_1,\ldots,\alpha_{t+t'})=\sum_{\{a_1,\ldots ,a_d\}\subset\F _p}\prod_{j=1}^{t+t'}\chi^{\alpha_j}
\big (\prod_{m=1}^d(b_{i_j}-a_m)\big ).$$
Nous appliquons maintenant le lemme \ref{somcar} pour majorer les sommes sur chaque $a_d$.
Cependant, le fait que les $a_j$ soient deux à  deux distincts rend les calculs un peu plus difficiles.

Dans une première étape, on écrit
\begin{equation}
\label{step1}
\begin{split}
W(\alpha_1,\ldots,\alpha_{t+t'})&=\sum_{\{a_1,\ldots ,a_{d-1}\}\subset\F _p}\prod_{j=1}^{t+t'}\chi^{\alpha_j}
\big (\prod_{m=1}^{d-1}(b_{i_j}-a_m)\big )\\
&\times \sum_{a_d\in\F_p\setminus\{ a_1,\ldots,a_{d-1}\}}
\chi ^{\alpha _j}\Big ( \prod_{ j=1}^{t+t'}(b_{i_j}-a_d)\Big )\\
& = \sum_{\{a_1,\ldots ,a_{d-1}\}\subset\F _p}\prod_{j=1}^{t+t'}\chi^{\alpha_j}
\big (\prod_{m=1}^{d-1}(b_{i_j}-a_m)\big )\sum_{a_d\in\F_p}
\chi ^{\alpha_j}\Big ( \prod_{ j=1}^{t+t'}(b_{i_j}-a_d)\Big )\\
&-\sum_{h=1}^{d-1}\sum_{\{a_1,\ldots ,a_{d-1}\}\subset\F _p}\prod_{j=1}^{t+t'}\chi^{\alpha_j}
\big (\prod_{\substack{m=1\\m\not =h}}^{d-1}(b_{i_j}-a_m)\big )\chi ^{2\alpha _h}\big  (\prod_{m=1}^{d-1}(b_{i_h}-a_m)\big )
\end{split}
\end{equation}
Dans cette deuxième ligne la variable $a_d$ est indépendante des autres.
En itérant ce procédé au bout de $d$ étapes, on obtient un nombre fini (au plus $d!$) de sommes de la forme
\begin{equation}
\label{Tal}
T(\boldsymbol{\alpha , \lambda}):=\Big (\sum_{a_1\in\F_p}\prod_{j=1}^{t+t'}\chi (b_{i_j}-a_1)^{\alpha_j\lambda_1}\Big )\cdots\Big (
\sum_{a_h\in\F_p}\prod_{j=1}^{t+t'}\chi (b_{i_j}-a_h)^{\alpha_j\lambda_h}\Big ),
\end{equation}
où $1\le h\le d$ et $\lambda_1,\ldots ,\lambda_h$ sont des entiers strictement positifs tels que
$$\lambda_1+\cdots +\lambda _h=d.$$
Dans la suite nous utiliserons le fait
que
\begin{equation}
\label{relh}
h=d-\sum_{m=1}^h(\lambda _m-1).
\end{equation}

Notons le pgcd $\alpha =(\alpha _1,\ldots ,\alpha_{t+t'})$ et $\alpha'_i=\alpha_i/\alpha$.
Pour chaque $m\in \{ 1,\ldots ,h\}$,
  on a

$$\sum_{a_m\in\F_p}\prod_{j=1}^{t+t'}\chi (b_{i_j}-a_m)^{\alpha_j\lambda_m}
=\sum_{a_m\in\F_p}\chi ^{\alpha\lambda_m}\Big ( \prod_{j=1}^{t+t'}(b_{i_j}-a_m)^{\alpha'_i}\Big ).$$
On peut appliquer le Lemme \ref{somcar} si $\chi^{\alpha\lambda_m}\not =\chi_0$.
Comme $(\alpha_1',\ldots ,\alpha'_{t+t'})=1$,
on a 
\begin{displaymath}\sum_{a_m\in\F_p}\prod_{j=1}^{t+t'}\chi (b_{i_j}-a_m)^{\alpha_j\lambda_m}
\le\left\{  \begin{array}{ll} (k-1)\sqrt{p}  & \text{si} \ \alpha\lambda_m\not\equiv 0\mod (p-1)\\
p & sinon.
\end{array}\right.
\end{displaymath}

Comme $(p-1)\not |\alpha$, lorsque $\lambda _m=1$, la somme est de module inférieur à 
$(k-1)\sqrt{p}$. Dans le cas où $\lambda_m\ge 2$, nous majorons la somme trivialement par $p$.
Notons $u(\boldsymbol{\lambda })$ le nombre d'indices $m$ tels que $\lambda_m\ge 2$.
On a
$$T(\boldsymbol{\alpha , \lambda})\le ((k-1)\sqrt{p} )^{h-u(\boldsymbol{\lambda })}
p^{u(\boldsymbol{\lambda })}\le (k-1)^d\sqrt{p}^{(h+u(\boldsymbol{\lambda }))}.$$
Mais d'après \eqref{relh}, $h+u(\boldsymbol{\lambda })\le d$
puisque $u(\boldsymbol{\lambda })\le \sum_{1\le m\le h}(\lambda _m -1)$.
Cela prouve que $T(\boldsymbol{\alpha,\lambda})\le (k-1)^dp^{d/2}$.

On obtient alors
\begin{equation}
\begin{split}
M&:=\Big | \varphi (p)^kT_2 (S,d)-{p-k\choose d}|S\setminus\{ 0\}|^\ell |S^c\setminus\{0\}|^{k-\ell}\Big | \\
&\ll
d!(k-1)^dp^{d/2}
\sum_{\substack {0\le t\le \ell\\ 0\le t'\le k-\ell\\(t ,t')\not =(0,0)}}
\sum_{ \substack{1\le i_1<\cdots <i_t\le\ell\\ \ell<i_{t+1}<\cdots<i_{t+t'}\le k}}
|S\setminus\{ 0\}|^{\ell -t}|S^c\setminus\{ 0\}|^{k-\ell -t'}\\
&\times\sum_{\chi_1\not =\chi_0}\cdots\sum_{\chi_{t+t'}\not=\chi_0}\prod_{j=1}^t\Big | \sum_{s\in S} \ov \chi_j (s)\Big |
\prod_{j=t+1}^{t+t'}\Big | \sum_{s\in S^c}\ov\chi _j (s)\Big |.
\end{split}
\end{equation}
Lorsque $S$ et $S^c$ sont de taille suffisamment grande, on peut obtenir des majorations intéressantes des sommes de caractères sur $S$ et $S^c$.
On commence par appliquer l'inégalité de Cauchy-Schwarz :
\begin{equation}
\begin{split}
\sum_{\chi\not =\chi_0}\big | \sum_{s\in S}\chi (s)\big | &=
 \sum_{\chi}\big | \sum_{s\in S}\chi (s)\big |-|S\setminus\{ 0\}\big |\\
&\le p^{1/2}\Big (\sum_\chi \big |\sum_{s\in S}\chi (s)\big |^2\Big )^{1/2}
-|S\setminus\{ 0\} |.
\end{split}
\end{equation}
On développe le carré puis on profite de l'orthogonalité des caractères :
$$
\sum_\chi \big |\sum_{s\in S}\chi (s)\big |^2  =\sum_\chi
\sum_{s_1,s_2\in S}\chi (s_1)\ov\chi (s_2)=\varphi (p)|S\setminus\{ 0\}|.$$
On en déduit l'inégalité :
$$\sum_{\chi\not =\chi_0}\big |\sum_{s\in S}\chi (s)|\le
\sqrt{p\varphi (p)|S\setminus\{ 0\}}\le p\sqrt{|S\setminus\{ 0\}|}.$$
On obtient de la même fa\c con une majoration de $\sum_{\chi\not =\chi_0}\Big | \sum_{s\in S^c}\chi (s)\Big |$.
Finalement,
\begin{equation}
\begin{split}
M&\ll d!(k-1)^dp^{d/2}\sum_{\substack {0\le t\le \ell\\ 0\le t'\le k-\ell\\
(t,t')\not =(0,0)}}
|S|^{\ell -t}|S^c|^{k-\ell -t'}
(p\sqrt{|S|})^t(p\sqrt{|S^c|})^{t'}
\sum_{ \substack{1\le i_1<\cdots <i_t\le\ell\\ \ell<i_{t+1}<\cdots<i_{t+t'}\le k}}1\\
&\ll d!(k-1)^dp^{d/2}\big [(|S|+p\sqrt{|S|})^\ell (|S^c|+p\sqrt{|S^c|})^{k-\ell}\\
&\ll d!(k-1)^d2^kp^{k+d/2}|S|^{\ell /2}|S^c|^{(k-\ell )/2}.
\end{split}
\end{equation}
Remarquons que le terme principal est de l'ordre de $\approx p^d|S^\ell|
|S^c|^{k-\ell}$.
Quitte à  échanger les rôles de $S$ et $S^c$, on peut supposer que
$|S|\le |S^c|$.

Dans ce cas la majoration que nous venons d'obtenir est la plus mauvaise quand
$\ell =k$.
Notre majoration est alors pertinente si
$$p^k\ll_{k,\ell ,d} p^{d/2}|S|^{k/2},$$
où la constante implicite dépend de $k,\ell ,d$.
Lorsque $|S|=O(1)$, cela impose  $k<d/2$.
Lorsque $|S|\gg p$, cela impose $k<d$.
Cela termine la preuve du théorème.

\section{La complexité $K_2$, preuve du Théorème \ref{k2general}}

Pour tout sous-ensemble $\cA\subset\F_p$ et tout $x\in\F_p$ nous utiliserons la notation suivante~:
$$x+\cA =\{ x+a:a\in\cA\}=\cA +x\quad x\cA =\{ xa: a\in\cA\}=\cA x.$$

Soit $\cA \subset \F _p$ de la forme
$\cA=\cB\cup\cC$ avec $\cB =\{ a_1,\ldots , a_\ell\}$ et
$\cC =\{ a_{\ell +1},\ldots ,a_{d+2}\}$, où $0\le\ell\le d+2$, de sorte que $|\cA|=d+2$.
On veut trouver $g\in\F_p [X]$ de degré plus petit ou égal à  $d$, sans racine multiple, tel que $g(a_i)\in S$ pour tout $1\le i\le \ell$
et $g(a_i)\in S^c$ pour tout $\ell +1\le i\le d+2$. Nous allons considérer le cas $\ell=1$ puis le cas $\ell\geq 2$.

\medskip

$\bullet$ Supposons d'abord $\ell=1$. Notons pour chaque $i\in\{2,3,\ldots,d+2\}$ le polynôme d'interpolation $f_i$ de degré plus petit ou égal à  $d$ tel que
\begin{equation*}
\left\{ \begin{array}{ll}  f_i (a_1)=0 \\
f_i (a_j)=1  \quad\text{si} \quad 2\le j\le d+2,\,j\neq i.
\end{array}\right.
\end{equation*}

Si $f_{i}(a_{i})= 0$ pour tout $i\in \{2,3,\ldots,d+2\}$, on forme alors le polynôme $f=f_2+f_3+\cdots+f_{d+2}$. On a
\begin{equation*}
\left\{ \begin{array}{ll}  f(a_1)=0 \\
f(a_i)=d+1\neq 0
 \quad\text{si} \quad i\in \{2,3,\ldots,d+2\}.
\end{array}\right.
\end{equation*}
On fixe alors $v\in S\setminus\{0\}$ et on choisit $u\not\in R(v)$, où
$$R(v)=\{ -v\ov{f(z_j)}, j \ \text{tel que}  f(z_j)\not =0\}\cup
\{(s-v)\ov{(d+1)} : s\in S\},$$
$z_1,\ldots ,z_m$, $m\le d-1$, représentant les racines distinctes du polynôme $f'(x)$. Cela est possible car $|R(v)|\le d-1+|S| <p$. Le polynôme $g(x)=uf(x)+v$ répond alors au problème.

Sinon $f_{i_0}(a_{i_0})\neq 0$ pour un certain $i_0\in \{2,3,\ldots,d+2\}$; on choisit alors $g$ de la forme $g(x)=uf_{i_0}(x)+v$ en fixant $v\in S\setminus\{0\}$ et en choisissant $u\not\in R(v)$,
$$R(v)=\{ -v\ov{f_{i_0}(z_j)}, j \ \text{tel que}  f_{i_0}(z_j)\not =0\}\cup
\{s-v : s\in S\}\cup\{(s-v)\ov{f_{i_0}(a_{i_0})} : s\in S\},$$
où $z_1,\ldots ,z_m$, $m\le d-1$, représentent les racines distinctes du polynôme $f'(x)$. Cela est possible car
$|R(v)|\le d-1+2|S| <p$.

\medskip

$\bullet$ Considérons dorénavant $\ell\geq 2$. Notons pour tout $i\in\{1,2,\ldots,\ell\}$ le polynôme d'interpolation $f_i$ de degré plus petit ou égal à  $d$ tel que
\begin{equation*}
\left\{ \begin{array}{ll}  f_i (a_j)=0 & \text{si}\quad 1\le j\le\ell,\,j\neq i\\
f_i (a_j)=1 & \text{si} \quad\ell+1\le j\le d+2.
\end{array}\right.
\end{equation*}

Si $f_{i}(a_{i})= 1$ pour tout $i\in \{1,2,\ldots,\ell\}$, alors le polynôme $f=f_1+f_2+\cdots+f_{\ell}-1$ vérifie
\begin{equation*}
\left\{ \begin{array}{ll}  f(a_j)=0 & \text{si}\quad 1\le j\le\ell\\
f(a_j)=\ell-1\neq0 & \text{si} \quad\ell+1\le j\le d+2,
\end{array}\right.
\end{equation*}
et l'on conclut comme dans un cas précédent en prenant $v\in S\setminus\{0\}$ et en choisissant $u\not\in R(v)$, où
$$R(v)=\{ -v\ov{f(z_j)}, j \ \text{tel que}  f(z_j)\not =0\}\cup
\{(s-v)\ov{(\ell-1)} : s\in S\}.$$ Le polynôme $g(x)=uf(x)+v$ convient alors.

Sinon, il existe $i_0\in \{1,2,\ldots,\ell\}$ tel que $\alpha:=f_{i_0}(a_{i_0})\neq 1$. Ainsi on désire trouver $g(x)=uf_{i_0}(x)+v$ (qui est de degré inférieur ou égal à  $d$) sans racine multiple vérifiant
\begin{equation*}
\left\{ \begin{array}{lll}  v\in S\\
\alpha u+v \in S\\
u+v \notin S\,.
\end{array}\right.
\end{equation*}
On peut supposer $\alpha\neq 0$ car autrement les deux premières équations sont équivalentes et on a déjà  vu qu'un tel polynôme $g$ existait. Notons $s=v$ et $s'=\alpha u+v$. Alors
\begin{equation*}
\left\{ \begin{array}{lll}  v\in S\\
\alpha u+v \in S\\
u+v \notin S
\end{array}\right.
\Longleftrightarrow
\left\{ \begin{array}{ll}  s,s'\in S\\
(1-\overline{\alpha})s+\overline{\alpha}s' \notin S\;.
\end{array}\right.
\end{equation*}
$g$ sans racine multiple signifie que $uf_{i_0}(z_j)+v=(1-\overline{\alpha}f_{i_0}(z_j))s+\overline{\alpha}f_{i_0}(z_j)s'\neq 0$ pour les racines $z_j$ de $f'_{i_0}$. Il n'est pas difficile de voir que l'égalité $(1-\overline{\alpha}f_{i_0}(z_j))s+\overline{\alpha}f_{i_0}(z_j)s'=0$ est vérifiée pour au plus $|S|$ couples $(s,s')$. Ainsi une condition suffisante pour terminer la preuve du théorème est~:
\begin{equation}\label{s1s2}
|\{ (s,s')\in S^2: (1-\overline{\alpha})s+\overline{\alpha}s' \notin S\}|> |S|(d-1)\,.
\end{equation}
Soit les ensembles $S_1=(1-\overline{\alpha})S$ et $S_2=\overline{\alpha}S$ de cardinalité $|S|$ (puisque $\alpha$ est différent de $0$ et $1$). Pour $n\in\F _p$ notons $r(n)$ le nombre de représentations $n=s_1+s_2$
avec $s_1\in S_1$ and $s_2\in S_2$.
D'après un résultat de Green and Ruzsa (\cite{GR}, Proposition 6.1) qui est une généralisation d'un théorème de Pollard, on a, pour tout $t\leq |S|$
\begin{eqnarray*}
\sum_{n\in\F _p}\min (t,r(n)) & \ge t\min (p, |S_1|+|S_2|-1-t)=t\min(p,2|S|-1-t)\\
&\ge t(2|S|-1-t) \quad\text{car}\quad 2|S|-1-t\leq p.
\end{eqnarray*}
Par ailleurs
$$\sum_{n\in S}\min (t,r(n))\leq t |S|\,.$$
Donc
\begin{eqnarray*}
\sum_{n\notin S} r(n)& \geq &\sum_{n\notin S}\min (t,r(n))= \sum_{n\in \F _p}\min (t,r(n))-\sum_{n\in S}\min (t,r(n))\\& \ge & t(2|S|-1-t)-t|S|=t(|S|-1-t)=:\phi(t)\,.
\end{eqnarray*}
La minoration est optimale pour $t_0=\frac{|S|-1}{2}$. $t_0$ ou $t_0+1/2$ étant un entier, on obtient
$$\sum_{n\notin S} r(n) \geq \phi(t_0+1/2)=\frac{|S|}{2}\left(\frac{|S|}{2}-1\right)\;.$$
Or $\frac{|S|}{2}\left(\frac{|S|}{2}-1\right)>|S|(d-1)$ car par hypothèse $|S|>4d+2$. En définitive on a bien l'inégalité voulue \eqref{s1s2}.

\noindent {\bf R. Balasubramanian}

\noindent Institute of Mathematical Sciences, C. I. T. Campus, Taramani, Chennai 600113, India

\noindent {\it e-mail :} balu@imsc.res.in

\medskip
\noindent {\bf C\'ecile Dartyge}

\noindent Institut \'Elie Cartan, Université de Lorraine, BP 70239, 54506 Vand\oe uvre-lès-Nancy Cedex, France

\noindent {\it e-mail :} cecile.dartyge@univ-lorraine.fr
\medskip

\noindent {\bf \'Elie Mosaki}

  \noindent Université de Lyon, Université Lyon 1, Institut
Camille Jordan CNRS UMR 5208, 43, boulevard du 11 Novembre 1918, F-69622
Villeurbanne

\noindent  {\it e-mail :} mosaki@math.univ-lyon1.fr

\end{document}